\documentclass[12pt,a4paper]{article}

\usepackage[T1]{fontenc}

\usepackage[french,american]{babel}

\newcommand {\citefr}[1]{\foreignlanguage{french}{#1}}

\usepackage{amsfonts}
\newcommand {\sC} {\ensuremath {\mathbb C}}
\newcommand {\sD} {\ensuremath {\mathbb D}}
\newcommand {\sE} {\ensuremath {\mathbb E}}
\newcommand {\sF} {\ensuremath {\mathbb F}}
\newcommand {\sN} {\ensuremath {\mathbb N}}
\newcommand {\sR} {\ensuremath {\mathbb R}}
\newcommand {\sZ} {\ensuremath {\mathbb Z}}

\DeclareMathSymbol{\carrefp}{\mathord}{AMSa}{3}
\newcommand {\finpreuve}{{\hspace*{\fill}} $\scriptstyle{\carrefp}$}

\DeclareMathSymbol{\restriction}{\mathord}{AMSa}{22}

\renewcommand {\labelenumi}{(\roman{enumi})}

\newtheorem{definition}{Definition}[section]
\newtheorem{thm}[definition]{Theorem}
\newtheorem{thinv}{Inversion theorem}

\newtheorem{prop}[definition]{Proposition}
\newtheorem{lemme}[definition]{Lemma}
\newtheorem{exemple}{Example}

\newcommand {\ssi} {if and only if }

\newcommand {\demrapide} [1]{\textbf {Proof:}~{#1}\finpreuve}

\newenvironment{dem}{\textbf {Proof: } \begin {quote}}
	{\finpreuve \end {quote}}

\newcommand {\etoile} {{\displaystyle \star}}
\newcommand {\quotient}[2] {\frac{\textstyle {#1}} {\textstyle {#2}}}
\newcommand {\qexp}[2] {\frac{\scriptstyle {#1}}{\scriptstyle {#2}}}
\newcommand {\epssur}[1] {\frac{\textstyle{\epsilon}}{\textstyle {{#1}}}}
\newcommand {\undemi} {\frac{\textstyle {1}}{\textstyle {2}}}
\newcommand {\et} {\:{\mathrm{and}}\:}
\newcommand {\kdn} {{k\in\sN}}
\newcommand {\ndn} {{n\in\sN}}
\newcommand {\ninf} {{n\to\infty}}
\newcommand {\xinx} {{x\in\grx}}
\newcommand {\ifg}[2] {\left[{#1},{#2}\right[}
\newcommand {\io}[2] {\left]{#1},{#2}\right[}
\newcommand {\ifd}[2] {\left]{#1},{#2}\right]}
\newcommand {\iseg}[2] {\left[{#1},{#2}\right]}
\newcommand {\pt} {\;\mbox{\textrm{\textit{.}}}}
\newcommand {\vrg} {\;\mbox{\textrm{\textit{,}}}}
\newcommand {\implique} {\,\Rightarrow\,}
\newcommand {\logequiv} {\,\Leftrightarrow\,}

\newcommand {\rplus} {\ensuremath{\sR_+}}
\newcommand {\rpet} {\ensuremath {\sR^\etoile_+}}

\newcommand {\alin}[2] {\ensuremath {\mathcal{L}_0\left({#1},{#2}\right)}}
\newcommand {\alinco}[2] {\ensuremath {\mathcal{L}\left({#1},{#2}\right)}}

\newcommand {\enlin}[1] {\ensuremath {\mathcal{L}_0\left({#1}\right)}}
\newcommand {\enlinco}[1] {\ensuremath {\mathcal{L}\left({#1}\right)}}
\newcommand {\couple} [2] {\ensuremath {\left\langle{#1},{#2}\right\rangle}}
\newcommand {\fctc} [2] {\ensuremath {\mathcal{C}\left({#1},{#2}\right)}}
\newcommand {\im}[1] {\ensuremath {\mathrm{R}\ep{#1}}}
\newcommand {\lp}[1] {\ensuremath{{\mathcal{L}}^{#1}}}
\newcommand {\ep}[1] {\ensuremath {\left({#1}\right)}}
\newcommand {\eva}[1] {\ensuremath {\left|{#1}\right|}}
\newcommand {\enor}[1] {\ensuremath {\left\|{#1}\right\|}}
\newcommand {\plp} [1] {\ensuremath {\ell^{#1}}}
\newcommand {\enp} [2] {\ensuremath {\left\|{#1}\right\|_{#2}}}
\newcommand {\enppuip} [2] {\ensuremath {\left\|{#1}\right\|_{#2}^{#2}}}
\newcommand {\eninf}[1] {\ensuremath {\left\|{#1}\right\|_{\infty}}}

\newcommand {\resol}{{R_T}}
\newcommand {\fun} {\ensuremath{\mathbf{1}}}
\newcommand {\intint}[2] %
	{\ensuremath {\displaystyle{\int\nolimits^{{#2}}_{{#1}}}\,}}

\newcommand {\gA} {\ensuremath {\mathcal {A}}}
\newcommand {\gT} {\ensuremath {\mathcal {T}}}
\newcommand {\gU} {\ensuremath {\mathcal {U}}}
\newcommand {\grx} {\ensuremath{\mathbf {X}}}

\newfont {\gothic} {eufm10 scaled 1200}

\newcommand {\ghC} {\ensuremath {\mbox{\gothic C}}}
\newcommand {\ghD} {\ensuremath {\mbox {\gothic D}}}

\newcommand {\refenum}[1] {\ensuremath {\mbox {\textit {({#1})}}}}

\begin{document}

\title{Some odd spectra}

\author {Philippe Deleval}

\date{October 2011}

\maketitle

\begin{abstract}
This text do not provide new great ideas, but fills a somptuous gap
between the state of the art in Operator Theory and the available
documentation (\cite{dusc}) of PDE specialists (see
\cite{auscher}, \cite{mcintosh}).

The spectrum for an element of an unitary Banach Algebra is nonvoid, closed
and bounded (i.e. compact). This result applies naturally to
the Banach Algebra of bounded operators on a Banach Space.

In the theory of Partial Derivate Equations appear some
``unbouded operators'' which are in fact \emph{partial} operators
defined on a dense subspace of the involved Banach Space.
These operators can have some odd spectra. The property allowing
to build a resolvent and a spectrum is the closeness. 

We want to show the problem is not in the fact that these operators
are unbounded, as suggested by \cite{dusc} (VII \S\ 1 N. 9),
but only in the fact they are partial.

We give briefly the methods for building resolvent ans spectra
for such operators and explain why the classical proofs do not work
in this case. After that we give some counterexamples, void or
unbounded spectra.
\end{abstract}

\tableofcontents

\section*{Notations}
\markboth{SOME ODD SPECTRA}{NOTATIONS}

\begin{itemize}
\item I write \couple{x}{y} for the ordered pair built with $x$ and $y$,
because I hold that the $(x,y)$ notation is far too much employed.

On the same basis, I write $\io{a}{b}$, french notation, rather than
$(a,b)$, the open interval $\{x\in\sR/a<x<b\}\pt$ In the
same way, I write $\ifg{a}{b}=\{x\in\sR/a\le x<b\}$ and
$\ifd{a}{b}=\{x\in\sR/a<x\le b\}\pt$
\item if $f$ is a linear mapping from the vector space \sF\ to
another vector space  \sE, I write $\ker\,f$ for the kernel of $f$,
\im{f} for its range.
\item \sE\ and \sF\ being vector spaces, I write \alin{\sF}{\sE}
for the set of linear mappings from \sF\ to \sE. If $\sF=\sE$,
I write \enlin{\sE} for \alin{\sE}{\sE}. $I$ is written for
an identical map, context must be clear enough to choose which identity.
\item if \sE\ and \sF\ are endowed with vector space topologies,
the continuous linear mapping (which are called ``bounded''
\emph{in the case of Banach Spaces}) are by definition the elements
of \alinco{\sF}{\sE} --- \enlinco{\sE} if $\sF=\sE$.
\item Banach spaces norms are generally written as $x\,\mapsto\,\enor{x}$,
context ought to make clear which norm (i.e. which space) is in question.

This is the case, among others, for operator norms,
always defined, for $T\in\alinco{\sE}{\sF}$
by $\enor{T}=\sup\left\{\enor{Tx}/x\in\sE\et\enor{x}\le1\right\}\pt$
\end{itemize}

When $T\in\alin{\sF}{\sE}$,
$G(T)=\{\couple{u}{v}\in\sF\times\sE/Tu=v\}$ (graph of $T$).

\section{Introduction: classical formulas}
\label{sec-intro}

Let \gA\ a Banach Algebra with unit $e$. The series expansion, for $x\in\gA$
such that $\enor{x}<1$, $(e-x)^{-1}=\sum\limits_\ndn\,x^n$
leads to the following result:

\begin{lemme}
\label{lemme-inversion-alba}
Let $x\in\gA$ such that $\enor{x}<1\pt$

Then $e-x$ is invertible in \gA\ and:
\begin{eqnarray*}
\enor{\ep{e-x}^{-1}} & \le & \quotient{1}{1-\enor{x}}\\
\enor{\ep{e-x}^{-1}-e} & \le & \quotient{\enor{x}}{1-\enor{x}}\\
\enor{\ep{e-x}^{-1}-e-x} & \le & \quotient{\enor{x}^2}{1-\enor{x}}\\
\end{eqnarray*}
\end{lemme}

The proof uses the power expansion $(e-x)^{-1}=\sum\limits_\kdn\,x^k$,
hence, for $\ndn$, \begin{math}(e-x)^{-1}-\sum\limits^n_{k=0}\,x^k=%
\sum\limits_{k>n}\,x^k=x^{n+1}(e-x)^{-1}\pt\end{math}

Note that the third inequality is not quite classical, but
the proof is fully similar to the proof of the second inequality.

By applying lemma \ref{lemme-inversion-alba} to $xa^{-1}$, we obtain
the following inversion theorem:

\begin{thinv}
\label{thinv-base}
Let $a\in\gA$ invertible.

Let $x\in\gA$ such that $\enor{x}<\enor{a^{-1}}^{-1}\pt$

Then $a-x$ is invertible in \gA\ and:
\begin{eqnarray*}
\enor{\ep{a-x}^{-1}} & \le & %
		\quotient{\enor{a^{-1}}}{1-\enor{x}.\enor{a^{-1}}}\\
\enor{\ep{a-x}^{-1}-a^{-1}} & \le & %
		\quotient{\enor{a^{-1}}^2.\enor{x}}{1-\enor{x}.\enor{a^{-1}}}\\
\enor{\ep{a-x}^{-1}-a^{-1}-a^{-1}xa^{-1}} & \le & %
		\quotient{\enor{a^{-1}}^3.\enor{x}^2}%
		{1-\enor{x}.\enor{a^{-1}}}\\
\end{eqnarray*}
\end{thinv}

This inversion theorem is the startpoint of definition of resolvents
and spectra in an unitary Banach Algebra, among others in \enlinco{\sE}
for a complex Banach space \sE.

Another inversion theorem, still straightforwardly proved
with the help of lemma \ref{lemme-inversion-alba}, is less usual.
I have seen it only in \cite{taf}, a book that Laurent Schwartz
has extracted from his course at \citefr{{\og}l'X\fg}
(\citefr{Ecole Polytechnique}), theorem T.2, XIV, 7; 1, p. 178 and
in \cite{mcintosh}, prop. 3.1.14 p.35. The wordings of these two authors
give only the first of the three inequalities.

\begin{thinv}
\label{thinv-op-Banachs}
Let \sE\ and \sF\ two Banach Spaces (both real or both complex).
Let $S\in\alinco{\sF}{\sE}$ having a inverse $S^{-1}$
in \alinco{\sE}{\sF}.

Let $T\in\alinco{\sF}{\sE}$ shuch that $\enor{T}<\enor{S^{-1}}^{-1}\pt$

Then $S-T\in\alinco{\sF}{\sE}$ is invertible in \alinco{\sE}{\sF} and:
\begin{eqnarray*}
\enor{\ep{S-T}^{-1}} & \le & %
		\quotient{\enor{S^{-1}}}{1-\enor{T}.\enor{S^{-1}}}\\
\enor{\ep{S-T}^{-1}-S^{-1}} & \le & %
		\quotient{\enor{S^{-1}}^2.\enor{T}}{1-\enor{T}.\enor{S^{-1}}}\\
\enor{\ep{S-T}^{-1}-S^{-1}-S^{-1}TS^{-1}} & \le & %
		\quotient{\enor{S^{-1}}^3.\enor{T}^2}%
		{1-\enor{T}.\enor{S^{-1}}}\\
\end{eqnarray*}
\end{thinv}

The proof applies the lemma \ref{lemme-inversion-alba}
in \enlinco{\sF} to $S^{-1}T$.

Like for classical inversion in a Banach algebra, these inequalities
show that the set of invertible $S\in\alinco{\sF}{\sE}$ is an open set
$\gU(\sF,\sE)$ and the inverse mapping
from $\gU(\sF,\sE)$ to $\gU(\sE,\sF)$ is a differentiable homeomorphism,
Its différential at $S$ is $T\,\mapsto\,-S^{-1}TS^{-1}$ --- element
of $\alinco{\alinco{\sF}{\sE}}{\alinco{\sE}{\sF}}\:$!

\section{Resolution between two Banach spaces}
\label{sec-res-two-Banachs}

In this section, we fix \sE\ and \sF\ to be two \emph{complex} Banach spaces,
$J$ a \emph{one-one} continuous mapping from \sF\ in \sE, such that
$\enor{J}\le1$. We explicitly assume that $J$ is not onto --- i.e. $J$
is not an isomorphism. The notation $I$ is for the identitiy (unit)
of \sE, $I'$ for the one of \sF.

\emph{Remark:} we may identify \sF\ with a subspace of \sE,
$J$ being the canonical injection. Il will be done in
section \ref{sec-th-spec-opfer}. Until then, the ``$J$'' is maybe somewhat
heavy in the notations, but gives more reliability when computing
and is responsible for the failure of the traditional proofs
that the spectrum is nonvoid and bounded.

\begin{definition}
\label{def-resolv-spectre}
Let $T\in\alinco{\sF}{\sE}\pt$

$\zeta\in\sC$ is resolved if $\zeta J-T$ has an inverse
$U\in\alinco{\sE}{\sF}$. $U$ will be called inverse (or reciprocal)
of $\zeta J-T$ and denoted $\resol(\zeta)\pt$

The set of resolved elements in \sC\ will be denoted $\rho(T)$
and called resolvent set of $T$, $\sC\setminus\rho(T)$ is
the \emph{spectrum} of $T$, denoted $\sigma(T)\pt$

The mapping assigning to $\zeta\in\rho(T)$ the operator
$\ep{\zeta J-T}^{-1}\in\alinco{\sE}{\sF}$ will be denoted
$\resol\,:\,\zeta\,\mapsto\,\resol(\zeta)$ and
called resolvent operator (or in brief resolvent) of $T$.
\end{definition}

N.B.: if, for $\zeta\in\sC$, $\zeta J-T$ has an algebraic inverse $U$
(i.e. in \alin{\sF}{\sE}), then by Banach's closed graph theorem, $U$
is bounded, so $\zeta$ is resolved.

One can prove as usual that $\rho(T)$ is open in \sC\ --- so $\sigma(T)$
is closed --- and, when $\rho(T)\neq\emptyset$, the resolvent function
is continuous and holomorphic from $\rho(T)$ to \alinco{\sE}{\sF}.

if $J$ is a linear isomorphism from \sF\ to \sE, we can prove
that $\sigma(T)$ is nonvoid and bounded.
This property is lost when $J$ is not invertible.

in fact, to prove that $\sigma(T)$ is nonvoid
and bounded, we use the invertibility of $J$ to show
that $J-\zeta^{-1}T$ is invertible if $\zeta$ is great enough.
So these properties cannot be proved by the usual means
for $J$ not invertible.
The counterexamples in the sections \ref{sec-ex-opfer}
and \ref{sec-ex-shift} will prove they are false.

\section{Spectral theory of the closed operators}
\label{sec-th-spec-opfer}

In this section, \sE\ is a complex Banach space, \sF\ a subspace
(\emph{a priori} not closed) of \sE. $J$ is the canonical injection
from \sF\ into \sE.

Because \sF\ is endowed with the induced topology, an unbounded
operator is nothing else than an element $T$ of
$\alin{\sF}{\sE}\setminus\alinco{\sF}{\sE}\pt$ We will omit in fact
the restriction $T\notin\alinco{\sF}{\sE}$! We fix a $T\in\alin{\sF}{\sE}\pt$
Now, the announced renorming.

\begin{prop}
\label{prop-pour-def-norme-assoc}
For $p\in[1,+\infty]$, the mapping $N_p^T$:
\begin{eqnarray*}
u & \mapsto & \ep{\enor{u}^p+\enor{Tu}^p}^\qexp{1}{p}\mbox{~if~}p<+\infty\vrg\\
u & \mapsto & \max\ep{\enor{u},\enor{Tu}}\mbox{~if~}p=+\infty\vrg
\end{eqnarray*}
is a norm. Moreover:
\begin{enumerate}
\item the norms $N_p^T$ are equivalent each other.
\item $T$ is continuous from \sF\ endowed with one of the normes $N_p^T$
to \sE.
\item the identity mapping from \sF\ endowed with the topology $\gT_2$
defined by the equivalent norms $N_p^T$ onto \sF\ endowed with
the induced topology $\gT_1$ is continuous. In other words, $\gT_2$
is stronger than $\gT_1$.
\item if $T$ is continuous from \sF\ with the induced topology to \sE,
then the $N_p^T$ norms are equivalent to the induced norm.
\end{enumerate}
\end{prop}
\begin{dem}
The relations $N_p^T(\alpha u)=\eva{\alpha}N_p^T(u)$,
$N_p^T(u)=0\logequiv u=0$ and the triangular inéquality are straightforward,
with the help of the Minkowski inequality for $p\in\io{1}{+\infty}\pt$
Other immediate relations : $\enor{u}\le N_p^T(u)$ (hence
\refenum{iii}), $\enor{Tu}\le N_p^T(u)$ (hence \refenum{ii}).

For $a$ and $b$ in \rplus\ and $p\in\ifg{1}{+\infty}$,
$\ep{a^p+b^p}^\qexp{1}{p}\le a+b$ follows from
Minkowski in $\sR^2$, with the two vectors \couple{a}{0} and \couple{0}{b}.
So, for $p\in[1,+\infty]$ and $u\in\sF$,
$N_\infty^T(u)\le N_p^T(u)\le N_1^T(u)\le2N_\infty^T(u)$, which proves
all the norms $N_p^T$ are equivalent to $N_\infty^T\pt$
Hence \refenum{i}.

Finally, if $T$ is continuous as soon as \sF\ is endowed with the induced
norm, we have $C\in\rplus$ such that $\enor{Tu}\le C\enor{u}$
for every $u\in\sF\pt$ Then, for $u\in\sF$,
$\enor{u}\le N_1^T(u)\le(1+C)\enor{u}$. That proves \refenum{iv}.
\end{dem}

\begin{definition}
\label{def-norme-assoc}
For $p\in[1,+\infty]$, the norm $N_p^T$
defined in prop. \ref{prop-pour-def-norme-assoc} is called
norm associated to $T$ with exponent $p$.

The topology defined by the $N_P^T$ norms will be called
topology associated to $T$.
We shall write $\gT_1$ for the induced topology on \sF,
$\gT_2$ for the topology defined by the equivalent norms $N_p^T$
\end{definition}

\emph{Remark:} it would have been a few easier to use merely
$p=1$ or $p=+\infty$. But $p=2$ will be usefull with Hilbert spaces!
Hence this little expensive generalization.

From \sF\ with this new norm into \sE, $T$ is continuous
(prop. \ref{prop-pour-def-norme-assoc} \refenum{ii}), and we
have another fact:

\begin{prop}
\label{prop-continuite-resolvante}
Let $\zeta\in\sC$ such that $\zeta J-T$
is \emph{algebraically} invertible, $S$ its inverse. Then:
\begin{enumerate}
\item $S$ is continuous from \sE\ to \sF\ endowed with $\gT_2\pt$
\item $S$ is continuous from \sE\ to \sF\ endowed with $\gT_1\pt$
\end{enumerate}
\end{prop}
\begin{dem}
\refenum{i} is an immediate consequence of Banach's closed graph theorem
(cf. the remark after def. \ref{def-resolv-spectre}).

\refenum{ii} follows because $\gT_1$ is coarser than $\gT_2\pt$
\end{dem}

To apply the results of section \ref{sec-res-two-Banachs}, it suffices
that \sF\ with $N_1^T$ be a Banach space.

The following theorem \ref{th-base-opfer} is very near prop. 3.1.4,
p. 34 of \cite{mcintosh}, given without proof.

\begin{thm}
\label{th-base-opfer}
Let $T\in\alin{\sF}{\sE}$, $G(T)$ its graph.
The following assertions are equivalent:
\begin{enumerate}
\item $G(T)$ is closed in $\sF\times\sE$.
\item for every sequence $(u_n)_\ndn$ with values in \sF, if $(u_n)$
converges to some $u\in\sE$ and if the sequence $(Tu_n)_\ndn$ converges
to some $v\in\sE$, then $u\in\sF$ et $v=Tu$.
\item \sF\ endowed with the topology associated to $T$ is a Banach space
with the norm $N_p^T$ for any $p\in[1,+\infty]$.
\end{enumerate}
\end{thm}
\begin{dem}
If $G(T)$ is closed and if the sequence $(u_n)_\ndn$ verifies the assumptions
of \refenum{ii}, then the sequence $\ep{\couple{u_n}{Tu_n}}_\ndn$ is
in $G(T)$ and converges to \couple{u}{v}, hence, since $G(T)$ is closed,
$\couple{u}{v}\in G(T)\pt$ In other words, $u\in\sF$ and $v=Tu$.
So $\refenum{i}\implique\refenum{ii}$.
\pagebreak[3.999999]

Suppose \refenum{ii}. Let $(u_n)_\ndn$ a Cauchy sequence
in \sF\ for, say, the norm $N_1^T\pt$

$\enor{u_p-u_q}\le N_1^T\ep{u_p-u_q}$
and $\enor{Tu_p-Tu_q}\le N_1^T\ep{u_p-u_q}$ for $p$ and $q$ in \sN, so
the sequences $(u_n)_\ndn$ and $(Tu_n)_\ndn$ are Cauchy in \sE\ and
converge, the first to a $u\in\sE$, the second to a $v\in\sE$.
$u\in\sF$ and $v=Tu$ by \refenum{ii}.

Then let $\epsilon\in\rpet$, $\eta=\epssur{2}\pt$
We have $N_1\in\sN$ such that, for $\ndn$ verifying $n\ge N_1$,
$\enor{u_n-u}<\eta$ and $N_2\in\sN$ such that, for $\ndn$ verifying
$n\ge N_2$, $\enor{Tu_n-Tu}<\eta\pt$
If $n\ge\max\ep{N_1,N_2}$, then:

\centerline{$N_1^T\ep{u-u_n}=\enor{u-u_n}+\enor{Tu-Tu_n}<2\eta=\epsilon\pt$}

Which proves that $(u_n)_\ndn$ converges to $u$ in \sF\ with $N_1^T$.
So, since the norms $N_p^T$ for $p\in[1,+\infty]$ are all
uniformly equivalent, \sF\ with an associated norm is complete.
Hence $\refenum{ii}\implique\refenum{iii}$.

Suppose now \refenum{iii}: \sF\ endowed with $N_1^T$
is therefore a Banach space.
Let $\ep{\couple{u_n}{Tu_n}}_\ndn$ in $G(T)$ converging to \couple{u}{v}
\emph{in $\sE\times\sE$ endowed with product topology}. Among others,
the sequence is Cauchy for this topology. One of the norms defining it
is $\nu\,:\,\couple{u}{v}\,\mapsto\,\enor{u}+\enor{v}\pt$

if $\epsilon\in\rpet$, then we have $N\in\sN$ such that, for every $p$ and $q$
in \sN\ verifying $p>N$ and $q>N$,
$\nu\ep{\couple{u_p}{Tu_p}-\couple{u_q}{Tu_q}}<\epsilon\pt$
in other words, $\nu\ep{\couple{u_p-u_q}{Tu_p-Tu_q}}<\epsilon\pt$

But that is still equivalent to $\enor{u_p-u_q}+\enor{T\ep{u_p-u_q}}<\epsilon$,
i.e. $N_1^T\ep{u_p-u_q}<\epsilon\pt$
So the sequence $(u_n)_\ndn$ is Cauchy in \sF\ with $N_1^T$,
supposed complete. $(u_n)_\ndn$ haa there a limit $u'$.

Moreover, since $\enor{u'-u_n}\le N_1^T\ep{u'-u_n}$,
$u'=\lim\limits_\ninf\,u_n=u$. And, since
$\enor{Tu'-Tu_n}\le N_1^T\ep{u'-u_n}$,
$Tu=Tu'=\lim\limits_\ninf\,Tu_n$ in \sE, hence $Tu=v$.

Which proves that $\couple{u}{v}\in G(T)$. So $G(T)$ is closed.
\end{dem} 

\begin{definition}
\label{def-op-ferme}
$T\in\alin{\sF}{\sE}$ is a closed operator if it verifies the equivalent
conditions of th. \ref{th-base-opfer}.
\end{definition}

If $T$ is a closed operator, whe can define its spectrum, its resolvent
because \sF\ is a Banach space with the topology $\gT_2$ and $T$ is bounded
for this topology. Like the icing of the cake, the continuity
of the resolvent does not depend on endowing \sF\ with the induced norm
or the new norm $N_1^T$ (prop. \ref{prop-continuite-resolvante}).
Elementary computation shows that $\resol(\zeta)$ is holomorphic
for both topologies.

\section{Three correlated examples}
\label{sec-ex-opfer}

Some examples will shed the light on the fickleness of spectra
for partial operators. With $T$ closed unbounded operator, we may have
$\rho(T)=\emptyset$ (\cite{mcintosh} p. 36, at the beginning of
section 3.2). So I'm keen on clearing out the ground. Thse
examples are very near of examples given in exercises
in  \cite{dusc}. 
The examples will be handled like propositions, because their assertions
have to be proved. Some computations preparing them come before.

In this section, $\grx=[0,1]$, $\ghC=\fctc{\grx}{\sC}$,
the set of continuous fonctions from \grx\ to \sC, a well-known \sC-based
separable Banach space. \ghD\ is the vector subspace of \ghC\ whose
elements are the continuously derivable mappings on \grx\ (right-derivable
at $0$, left-derivable at $1$). \ghC\ is endowed with the usual norm,
$f\,\mapsto\,\max\limits_\xinx\,\eva{f(x)}$, restriction to \ghC\ of
the \lp{\infty}-seminorm. It will be denoted $f\,\mapsto\,\eninf{f}\pt$
$T$ will  be, for these three examples, the derivation operator
form \ghD\ to \ghC, or its restriction to a subspace \sF\ of \ghD.

\begin{lemme}
\label{lemme-calcul-aux-pour-ex}
Let $k\in\sC^\etoile$, $a$ and $b$ in $\sC$.

The mapping
$t\,\mapsto\,\ep{\quotient{at+b}{k}-\quotient{a}{k^2}}e^{kt}$
is a primitive of the mapping  from \sR\ to
\sC\ $t\,\mapsto\,\ep{at+b}e^{kt}\pt$
\end{lemme}
\demrapide{elementary calculus.}

\begin{definition}
\label{def-fcts-de-zeta}
Let $\zeta\in\sC$ A mapping and a mapping family, all these in \ghD,
are associated to $\zeta$:
\begin{enumerate}
\item $h_\zeta\,:\,x\,\mapsto\,e^{\zeta x}\pt$
\item for each $f\in\ghC$,
$K_{\zeta}f\,:\,x\,\mapsto\,e^{\zeta x}\intint{0}{x}e^{-\zeta t}f(t)dt\pt$
\end{enumerate}
\end{definition}


\begin{lemme}
\label{lemme-equa-commune-ex}
Let $\zeta\in\sC$, $f\in\ghC$, $a=\Re\zeta\pt$
{
\renewcommand {\labelenumii}{(\roman{enumii})}
\renewcommand {\labelenumi}{(\arabic{enumi})}
\begin{enumerate}
\item if $a\neq0$, $\enor{K_\zeta}=\quotient{e^a-1}{a}>0\pt$
\item if $a=0$ $\enor{K_\zeta}=1\pt$
\item For every $u\in\ghD$, the following assertions are equivalent:
\begin{enumerate}
\item $\zeta u-u'=f\pt$
\item there exists a constant $\gamma\in\sC$ such that
$u=\gamma h_\zeta-K_{\zeta}f\pt$
\end{enumerate}
\end{enumerate}}
\end{lemme}
\begin{dem}
$\eva{K_{\zeta}f(x)}=\eva{e^{\zeta x}\intint{0}{x}e^{-\zeta t}f(t)dt}\vrg$
so:

\centerline{\begin{math}\eva{K_{\zeta}f(x)}\le\eva{e^{\zeta x}}\intint{0}{x}%
\eva{e^{-\zeta t}}\eva{f(t)}dt\le%
e^{ax}\eninf{f}\intint{0}{x}e^{-at}dt\pt\end{math}}
\pagebreak[3.999999]

If $a\neq0$, the mapping $t\,\mapsto\,e^{-at}$ is the derivate
of $t\,\mapsto\,\quotient{e^{-at}}{-a}$, so
$\intint{0}{x}e^{-at}dt=\quotient{e^{-ax}-1}{-a}=\quotient{1-e^{-ax}}{a}\pt$
Notice that, for $a>0$, $1-e^{-ax}\ge0$ and, for $a<0$, $1-e^{-ax}\le0$:
the quotient is nonnegative.

It follows that \begin{math}\eva{K_{\zeta}f(x)}\le%
e^{ax}\quotient{1-e^{-ax}}{a}\eninf{f}=%
\quotient{e^{ax}-1}{a}\eninf{f}\pt\end{math}

If $a>0$, then $e^{ax}\le e^a$ for every $\xinx$,
so $\quotient{e^{ax}-1}{a}\le\quotient{e^a-1}{a}$, hence
$\eninf{K_{\zeta}f}\le\quotient{e^a-1}{a}\eninf{f}\pt$

If $a<0$, then $-e^{ax}\le -e^a$ for every $\xinx$ and consequently
\begin{math}\quotient{e^{ax}-1}{a}=\quotient{1-e^{ax}}{-a}\le%
\quotient{1-e^a}{-a}=\quotient{e^a-1}{a}\end{math}, hence, like in
the previous case, $\eninf{K_{\zeta}f}\le\quotient{e^a-1}{a}\eninf{f}\pt$

If $a=0$, $\intint{0}{x}e^{-at}dt=\intint{0}{x}1dt=x$ and $e^{ax}=1$,
hence straight $\eva{K_{\zeta}f(x)}\le x\eninf{f}$, with
implies $\eninf{K_{\zeta}f}\le\eninf{f}\pt$

So we have $\enor{K_\zeta}\le\quotient{e^a-1}{a}$ if $a\neq0$,
$\enor{K_\zeta}\le1$ if $a=0$.

Now, let $b=\Im\,\zeta$, so that $\zeta=a+bi$.
Let $f\in\ghC$ defined by $f(x)=e^{bix}\pt$
Since $\eva{f(x)}=1$ for every $\xinx$, $\eninf{f}=1$.
For $\xinx$:
\begin{eqnarray*}
K_{\zeta}f(x) & = & e^{\zeta x}\intint{0}{x}e^{-\zeta t}e^{bit}dt\\
	& = & e^{\zeta x}\intint{0}{x}e^{-at}dt\pt
\end{eqnarray*}
If $a=0$, $K_{\zeta}f(x)=xe^{bix}\vrg$ so $\eva{K_{\zeta}f(x)}=x$
and $\eninf{K_{\zeta}f}=1$. Which proves $\enor{K_\zeta}=1\pt$

If $a\neq0$:
\begin{eqnarray*}
K_{\zeta}f(x) & = & e^{ax}e^{bix}\quotient{e^{-ax}-1}{-a}\\
	& = & e^{ibx}\quotient{e^{ax}-1}{a}\\
\eva{K_{\zeta}f(x)} & = & \quotient{e^{ax}-1}{a}\pt
\end{eqnarray*}

Hence (see the above computations of upper bounds) we have
$\eninf{K_{\zeta}f}=\quotient{e^{a}-1}{a}$
and $\enor{K_\zeta}=\quotient{e^{a}-1}{a}\pt$

\refenum{3} is classical calculus.
\end{dem}

\begin{prop}
\label{prop-base-ex-syst-oper}
Let \sE\ a closed subspace of \ghC\ (therefore a Banach space with induced
norm and topology), $\sF=\{u\in\sE\cap\ghD/u'\in\sE\}\pt$
The operator $T$ from \sF\ to \sE\ is closed.
\end{prop}
\begin{dem}
Suppose $(u_n)_\ndn$ sequence in \sF\ converging in \sE\ to $u$,
suppose more that $u'_n$ converges in \sE\ to $v$.

A classical calculus result show that the sequence $(w_n)_\ndn$,
defined by $w_n(x)=\intint{0}{x}u'_n(t)dt\vrg$
converges in \ghC\ to $w(x)=\intint{0}{x}v(t)dt\pt$

Now, for $\xinx$, $\intint{0}{x}u'_n(t)dt=u_n(x)-u_n(0)\pt$
Since $(u_n)_\ndn$ converges uniformly to $u$,
$\lim\limits_\ninf\,\ep{u_n(x)-u_n(0)}=u(x)-u(0)\pt$

Therefore $u(x)-u(0)=\intint{0}{x}v(t)dt\vrg$ hence
$u\in\ghD$ and $u'=v$. Moreover, since $u'_n\in\sE$ for every $\ndn$
and \sE\ is closed in \ghC, $v=u'\in\sE$.
Which proves that the graph $G(T)$ of $T$ is closed, so $T$ is closed.
\end{dem}

\begin{prop}
\label{prop-T-noyau-flinc}
Let $\Lambda$ a continuous linear functional on \ghC. Let $\sE=\ker\,\Lambda$,
$\sF=\{u\in\sE\cap\ghD/u'\in\sE\}=\{u\in\ghD/\Lambda u=0\et\Lambda u'=0\}$
(cf. prop. \ref{prop-base-ex-syst-oper}).
Then \sE\ is closed, $T$ is a closed operator from \sF\ to \sE\ and,
for every $\zeta\in\sC$, the following assertions are equivalent:
\begin{enumerate}
\item $\zeta\in\sigma(T)\pt$
\item $h_\zeta\in\sE$, i.e. $\Lambda\ep{h_\zeta}=0\pt$
\end{enumerate}
Moreover, if $\zeta\notin\sigma(T)$, for every $f\in\sE$,
we have $\resol(\zeta)f=\gamma h_\zeta-K_{\zeta}f$ with
$\gamma=\quotient{\Lambda\ep{K_{\zeta}f}}{\Lambda\ep{h_\zeta}}\pt$
\end{prop}
\begin{dem}
\sE, kernel of the continuous linear functional $\Lambda$, is closed
and then prop. \ref{prop-base-ex-syst-oper} proves that $T$ is closed.

Let $\zeta\in\sC\pt$ We have, for given $f\in\sE$, to resolve
the equation in $u\in\sF$ $\zeta u'-u=f\pt$
By lemma \ref{lemme-equa-commune-ex}, for $u\in\ghD$, $\zeta u-u'=f$
\ssi there exists $\gamma\in\sC$ such that
$u(x)=\gamma e^{\zeta x}-K_{\zeta}f(x)\pt$
Can we have $u\in\sF$?

First, if $u\in\sE$, then $u\in\sF$ since, from $\zeta u-u'=f$,
we deduce $\Lambda\ep{u'}=\zeta\Lambda(u)-\Lambda(f)=0$
for $f\in\sE$ and $u\in\sE\cap\ghD$.

Let $A(\zeta)=\Lambda\ep{h_\zeta}$ and
$B(\zeta)=\Lambda\ep{K_{\zeta}f}\pt$

$\Lambda(u)=\gamma A(\zeta)-B(\zeta)\vrg$ so everything
depends on $A(\zeta)$:
\begin{itemize}
\item if $A(\zeta)\neq0$, there exists one and only one
$\gamma=\quotient{B(\zeta)}{A(\zeta)}$ such that $u\in\sE$, so $u\in\sF\pt$

And we have $M\in\rpet$ such that $\eva{\Lambda\ep{K_{\zeta}f}}\le M\eninf{f}$
by continuity of $\Lambda$ and lemma \ref{lemme-equa-commune-ex}, hence:
\begin{displaymath}
\eninf{u}\le\ep{1+\quotient{1}{\eva{A(\zeta)}}}M\eninf{f}\pt
\end{displaymath}
Which proves that $\zeta\in\rho(T)\pt$
\item if $A(\zeta)=0$, following the value of $B(\zeta)$, there is
no value of an infinity of values for $\gamma$ such that $u\in\sE$
(and $u\in\sF$). So $\zeta\in\sigma(T)\pt$
\end{itemize}

Hence $\sigma(T)=\{\zeta\in\sC/A(\zeta)=0\}\pt$
\end{dem}

Continuous linear fuctionals on \ghC\ are nothing else than functional Radon
measures on $[0,1]\pt$ We will write $\delta_x$ for the Dirac measure
on $x\in[0,1]$ ($\delta_x(f)=f(x)$ for $f\in\ghC$).

The three examples thereafter use prop. \ref{prop-T-noyau-flinc}.
Except for the first one (where $\rho(T)=\emptyset$), we consider
the behaviour when $\eva{\zeta}\to\infty$ of $\enor{\resol(\zeta)}\pt$

\begin{exemple}
\label{exemple1-opfer}
Let us take for $\Lambda$ the null functional.

$\sE=\ghC$ (hence $\sF=\ghD$), $T$ is closed, $\sigma(T)=\sC$
and $\rho(T)=\emptyset\pt$
\end{exemple}
\begin{dem}
By prop. \ref{prop-T-noyau-flinc}, $h_\zeta\in\ghC=\ker\,\Lambda$
for every $\zeta\in\sC$, so $\sigma(T)=\sC\pt$
\end{dem}

\emph{N.B:} in this example, $T$ is onto;
for $f\in\ghC$, $u\,:\,x\,\mapsto\,\intint{0}{x}f(t)dt$ is in \ghD\ and $Tu=f$.

\begin{exemple}
\label{exemple2-opfer}
Let us take $\Lambda=\delta_0\pt$

$\sE=\{f\in\ghC/f(0)=0\}$ is closed, $\sF=\{u\in\ghD/u(0)=0\et u'(0)=0\}$ and,
for every $\zeta\in\sC$, $\zeta J-T$ is a bijection with reciprocal
$\resol(\zeta)=-K_\zeta\vrg$ continuous mapping.

$\rho(T)=\sC$ and $\sigma(T)=\emptyset\pt$

For $\zeta\in\rpet$, \begin{math}\quotient{e^\zeta-1-\zeta}{\zeta^2}\le%
\enor{\resol(\zeta)}\le\quotient{e^\zeta-1}{\zeta}\pt\end{math}

So
\begin{math}\lim\limits_{{\zeta\to+\infty}\atop{\zeta\in\rplus}}\,%
\enor{\resol(\zeta)}=+\infty\end{math}
\end{exemple}
\begin{dem}
We apply prop. \ref{prop-T-noyau-flinc}.
Since $\Lambda\ep{h_\zeta}=1$ for each $\zeta\in\sC$,
$\sigma(T)=\emptyset\pt$

$K_{\zeta}f(0)=e^{\zeta x}\intint{0}{0}e^{-\zeta t}f(t)dt=0\vrg$
so $\Lambda\ep{K_{\zeta}f}=0$, and, using the formulas of
prop. \ref{prop-T-noyau-flinc}, we obtain $\gamma=0$ and
$\resol(\zeta)f=-K_{\zeta}f\pt$

So $\enor{\resol(\zeta)}=\enor{K_\zeta}$ (norms in \alinco{\sE}{\sF}).
We deduce, by lemma \ref{lemme-equa-commune-ex} (since $\sE\subseteq\sC$)
$\enor{K_\zeta}\le\quotient{e^\zeta-1}{\zeta}$ --- notice
that the operator norm of the derivation operator \emph{from \ghC}
has been computed by using the mapping $\fun\,:\,x\,\mapsto\,1$,
which is not in \sE.

To get a lesser bound of $\enor{K_\zeta}$, we will compute $\eninf{K_{\zeta}f}$
for some $f\in\sE$ such that $\eninf{f}=1$. We merely take for $f$
the identity of \grx, $f(x)=x$.

By lemma \ref{lemme-calcul-aux-pour-ex}, the mapping
$x\,\mapsto\,xe^{-\zeta x}$, which is involved in the computation,
admits the primitive $x\,\mapsto\,\ep{-\zeta^{-1}x-\zeta^{-2}}e^{-\zeta x}\pt$

For $\xinx$:
\begin{eqnarray*}
K_{\zeta}f(x) & = & e^{\zeta x}\intint{0}{x}te^{-zeta t}dt\\
	& = & e^{\zeta x}%
		\left[\ep{-\zeta^{-1}t-\zeta^{-2}}e^{-\zeta t}\right]^x_0\\
	& = & \quotient{-e^{\zeta x}}{\zeta^2}%
		\left[\ep{\zeta t+1}e^{-\zeta t}\right]^x_0\\
	& = & \quotient{-e^{\zeta x}}{\zeta^2}%
		\ep{\ep{\zeta x+1}e^{-\zeta x}-1}\\
	& = & \quotient{e^{\zeta x}-1-\zeta x}{\zeta^2}\pt
\end{eqnarray*}
The mapping $g\,:\,x\,\mapsto\,e^{\zeta x}-1-\zeta x$ is nondecreasing
(see it like partial summation of the exponential series,
or use the derivate,
$g'(x)=\zeta e^{\zeta x}-\zeta=\zeta\ep{e^{\zeta x}-1}\ge0$).

Since $g(0)=0$, $\eninf{g}=g(1)=e^\zeta-1-\zeta\pt$

Because $K_{\zeta}f=\quotient{1}{\zeta^2}g$:

\begin{displaymath}
\eninf{K_{\zeta}f}=\quotient{1}{\zeta^2}\eninf{g}=%
\quotient{e^\zeta-1-\zeta}{\zeta^2}\pt
\end{displaymath}

Since $\eninf{f}=1$, we obtain
$\enor{\resol(\zeta)}=\enor{K_\zeta}\ge\quotient{e^\zeta-1-\zeta}{\zeta^2}\pt$

It suffices now to observe that
\begin{math}\lim\limits_{{\zeta\to+\infty}\atop{\zeta\in\rplus}}\,%
\quotient{e^\zeta-1-\zeta}{\zeta^2}=+\infty\end{math}
\end{dem}

In the following example \ref{exemple3-opfer}, $\resol(\zeta)\neq -K_\zeta\pt$

\begin{exemple}
\label{exemple3-opfer}
Let $\Lambda$ the bounded linear functional
on \sC\ $\delta_\qexp{1}{2}-\delta_0$, such that, for $f\in\sC$,
$\Lambda f=f\ep{\undemi}-f(0)\vrg$
$\sE=\left\{f\in\ghC/f(0)=f\ep{\undemi}\right\}\vrg$ kernel of $\Lambda$.

\sE\ is closed and $\sigma(T)=4i\pi\sZ=\{4i\pi n/n\in\sZ\}\pt$

For $f\in\sE$, $\zeta\in\rho(T)$ and $\xinx$:
\begin{displaymath}
\resol(\zeta)f(x)=\quotient{e^{\zeta x}}{e^\qexp{\zeta}{2}-1}%
K_{\zeta}f\ep{\undemi}-K_{\zeta}f(x)\pt
\end{displaymath}

Let $\zeta\in\rpet$ (hence $\zeta\in\rho(T)$). Then there exists a mapping
$f\in\sE$ such that:
\begin{enumerate}
\item $f(0)=f\ep{\undemi}=0\pt$
\item $\eninf{f}=e^\qexp{\zeta}{2}\pt$
\item $K_{\zeta}f\ep{\undemi}=0\pt$
\item for every $\xinx$, $\resol(\zeta)f(x)=-K_{\zeta}f(x)\pt$
\item \begin{math}\resol(\zeta)f(1)=e^\qexp{\zeta}{2}%
\ep{\quotient{1}{\zeta}+\quotient{2}{\zeta^2}}%
-\quotient{2e^\zeta}{\zeta^2}\pt\end{math}
\end{enumerate}
We have therefore:
\begin{eqnarray*}
\enor{\resol(\zeta)} & \ge & \quotient{2e^\qexp{\zeta}{2}}{\zeta^2}%
		-\quotient{1}{\zeta}-\quotient{2}{\zeta^2}\vrg\\
\lim\limits_{{\zeta\to+\infty}\atop{\zeta>0}}\,%
\enor{\resol(\zeta)} & = & +\infty\pt
\end{eqnarray*}
\end{exemple}
\begin{dem}
Since \sE\ is the kernel of $\Lambda$, prop. \ref{prop-T-noyau-flinc}
applies.
\pagebreak[3.999999]

Let $\zeta\in\sC\pt$
$\Lambda\ep{h_\zeta}=e^\qexp{\zeta}{2}-1\pt$
So $\zeta\in\sigma(T)$ \ssi $e^\qexp{\zeta}{2}=1$, i.e. \ssi there exists
$n\in\sZ$ such that $\quotient{\zeta}{2}=2i\pi n$,
i.e. $\zeta=4i\pi n\pt$

To compute $\gamma$ by using prop. \ref{prop-T-noyau-flinc} formula,
for $f\in\sE$, since $K_{\zeta}f(0)=0\vrg$ \begin{math}\Lambda\ep{K_{\zeta}f}=%
K_{\zeta}f\ep{\undemi}-K_{\zeta}f(0)=%
K_{\zeta}f\ep{\undemi}\pt\end{math}

Hence, for $\zeta\in\rho(T)=\sC\setminus4i\pi\sZ$, $\xinx$:
\begin{eqnarray*}
\gamma & = & \quotient{K_{\zeta}f\ep{\undemi}}%
		{e^\qexp{\zeta}{2}-1}\vrg\\
\resol(\zeta)f(x) & = & \quotient{e^{\zeta x}}{e^\qexp{\zeta}{2}-1}%
K_{\zeta}f\ep{\undemi}-K_{\zeta}f(x)\pt
\end{eqnarray*}

Let $\zeta\in\rpet$ fixed and $f$ the mapping defined
from $\grx=[0,1]$ in \sR\ by:
\begin{itemize}
\item if $x\in\iseg{0}{\undemi}$, $f(x)=e^{\zeta x}\sin\,4\pi x\pt$
\item if $x\in\ifd{\undemi}{1}$, $f(x)=e^\qexp{\zeta}{2}(2x-1)\pt$
\end{itemize}

So $f$ takes the  value of $f_1(x)=e^{\zeta x}\sin\,4\pi x$
for $x\le\undemi$ and the value of $f_2(x)=e^\qexp{\zeta}{2}(2x-1)$
for $x>\undemi\pt$ $f_1$ and $f_2$ are real analytic,
thus continuous. To establish the continuity of $f$, it suffices
to show $f_1\ep{\undemi}=f_2\ep{\undemi}\pt$ Now:
\begin{eqnarray*}
f_1\ep{\undemi} & = & e^\qexp{\zeta}{2}\sin\,2\pi\\
	& = & 0\pt\\
f_2\ep{\undemi} & = & e^\qexp{\zeta}{2}\ep{2\undemi-1}\\
	& = & 0\pt
\end{eqnarray*}
Hence moreover $f\ep{\undemi}=0$ and $f(x)=f_2(x)$
on $\iseg{\undemi}{1}\pt$

To fulfil \refenum{i}, it remains to compute
$f(0)=\sin\,4\pi0=\sin\,0=0\pt$

By the way, $f\ep{\undemi}=f(0)$, so $f\in\sE$.

For $x\in\iseg{0}{\undemi}$, \begin{math}\eva{f(x)}=%
e^{\zeta x}\eva{\sin\,4\pi x}\le e^{\zeta x}\le e^\qexp{\zeta}{2}\pt\end{math}

For $x\in\iseg{\undemi}{1}$, $0\le 2x-1\le1$, so
$0\le f(x)\le e^\qexp{\zeta}{2}\pt$
Thus, $\eninf{f}\le e^\qexp{\zeta}{2}\pt$

But $f(1)=e^\qexp{\zeta}{2}$, so $\eninf{f}=e^\qexp{\zeta}{2}\pt$

We compute $I(x)=\intint{0}{x}e^{-\zeta t}f(t)dt$ and
$K_{\zeta}f(x)=e^{\zeta x}I(x)$ in two steps:
\begin{itemize}
\item if $x\in\iseg{0}{\undemi}$:
\begin{eqnarray*}
I(x) & = & \intint{0}{x}e^{-\zeta t}e^{\zeta t}\sin\,4\pi tdt\\
	& = & \intint{0}{x}\sin\,4\pi tdt\\
	& = & \left[\quotient{-\cos\,4\pi t}{4\pi}\right]_0^x\\
	& = & \quotient{1-\cos\,4\pi x}{4\pi}\pt\\
K_{\zeta}f(x) & = & e^{\zeta x}\,\quotient{1-\cos\,4\pi x}{4\pi}\pt
\end{eqnarray*}

Observe by the way that $I\ep{\undemi}=\quotient{1-\cos\,2\pi}{4\pi}=0$
and thus $K_{\zeta}f\ep{\undemi}=0\pt$
\item if $x\in\iseg{\undemi}{1}$:
\begin{eqnarray*}
I(x) & = & \intint{0}{\qexp{1}{2}}e^{-\zeta t}f(t)dt+%
		\intint{\qexp{1}{2}}{x}e^{-\zeta t}f(t)dt\\
	& = & I\ep{\undemi}+\intint{\qexp{1}{2}}{x}e^{-\zeta t}f(t)dt\\
	& = & \intint{\qexp{1}{2}}{x}e^{-\zeta t}f(t)dt\\
	& = & e^\qexp{\zeta}{2}\,%
		\intint{\qexp{1}{2}}{x}(2t-1)e^{-\zeta t}f(t)dt\\
	& = & e^\qexp{\zeta}{2}\,%
		\left[\ep{\quotient{2t-1}{-\zeta}-\quotient{2}{\zeta^2}}%
		e^{-\zeta t}\right]^x_\qexp{1}{2}\\
	& = & e^\qexp{\zeta}{2}\ep{\ep{\quotient{1-2x}{-\zeta}-%
		\quotient{2}{\zeta^2}}e^{-\zeta x}+%
		\quotient{2e^\qexp{-\zeta}{2}}{\zeta^2}}\pt\\
K_{\zeta}f(x) & = & e^\qexp{\zeta}{2}\ep{\quotient{1-2x}{\zeta}-%
		\quotient{2}{\zeta^2}}+\quotient{2e^{\zeta x}}{\zeta^2}\pt\\
\end{eqnarray*}

Among other values, \begin{math}K_{\zeta}f(1)=e^\qexp{\zeta}{2}%
\ep{-\quotient{1}{\zeta}-\quotient{2}{\zeta^2}}+%
\quotient{2e^\zeta}{\zeta^2}\end{math} 
\end{itemize}

Since $K_{\zeta}f\ep{\undemi}=0$, the formula for computing
the resolvent becomes $\resol(\zeta)f(x)=-K_{\zeta}f(x)$, so:
\begin{itemize}
\item if $x\in\iseg{0}{\undemi}$,
$\resol(\zeta)f(x)=-e^{\zeta x}\,\quotient{1-\cos\,4\pi x}{4\pi}\pt$
\item if $x\in\iseg{\undemi}{1}$,
\begin{math}\resol(\zeta)f(x)=e^\qexp{\zeta}{2}\ep{\quotient{2x-1}{\zeta}+%
\quotient{2}{\zeta^2}}-\quotient{2e^{\zeta x}}{\zeta^2}\pt\end{math}
\end{itemize}

We deduce \begin{math}\resol(\zeta)f(1)=e^\qexp{\zeta}{2}%
\ep{\quotient{1}{\zeta}+\quotient{2}{\zeta^2}}%
-\quotient{2e^\zeta}{\zeta^2}\pt\end{math}
What remains to prove is deduced from:
\begin{eqnarray*}
-\resol(\zeta)f(1) & \le & \eva{\resol(\zeta)f(1)}\le\eninf{\resol(\zeta)f}\\
\enor{\resol(\zeta)} & \ge & \quotient{\eninf{\resol(\zeta)f}}{\eninf{f}}%
=\quotient{\eninf{\resol(\zeta)f}}{e^\qexp{\zeta}{2}}\pt
\end{eqnarray*}
\end{dem}

\section{Example with the shift}
\label{sec-ex-shift}

In the primitive text wherefrom this paper is derived, the fourth example was
again with a subspace \sE\ of $\ghC=\fctc{[0,1]}{\sC}$. But my impression
was that this example was basically a shift.

So, in this section, we fix a $p\in[1,+\infty]$ and work with the complex
Banach space $\sE=\plp{p}\pt$ An element $x$ of \plp{p} will be systematically
written $x=(x_n)_\ndn\pt$

\begin{definition}
\label{def-aux-ex-opcont}
$S$ is the mapping from \sE\ to \sE\ such that, for every $x\in\ghC$,
$Sx=y$ with $y_n=x_{n+1}\pt$
\end{definition}
\pagebreak[3.999999]

\begin{prop}
\label{prop-aux-ex-opcont}
$S$ from def. \ref{def-aux-ex-opcont} has the following properties:
\begin{enumerate}
\item $S\in\enlinco{\plp{p}}$ and $\enp{S}{p}=1\pt$
\item $S$ is onto.
\item $\sigma\ep{S}=\sD$, with $\sD=\{\zeta\in\sC/\eva{\zeta}\le1\}$
(unitary disk of \sC).
\item for $\zeta\in\rho\ep{S}$ (i.e. $\zeta\in\sC$ such that $\eva{\zeta}>1$),
$\resol(\zeta)(x)=y$, with $y$ such that, for every $\ndn$:
\begin{displaymath}
y_n=\sum\limits^{\infty}_{k=0}\,\zeta^{-k-1}x_{n+k}\pt
\end{displaymath}
Among others, $y\in\plp{p}\pt$
\end{enumerate}
\end{prop}
\begin{dem}
The linearity of $S$ is obvious.

Let $x\in\plp{p}$
\begin{itemize}
\item if $p<+\infty$, \begin{math}\enppuip{Sx}{p}=%
\sum\limits_\ndn\,\eva{x_{n+1}}^p%
=\sum\limits^{\infty}_{n=1}\,\eva{x_n}^p\pt\end{math}

Since $\enppuip{Sx}{p}+\eva{x_0}^p=\enppuip{x}{p}$,
$\enp{Sx}{p}\le\enp{x}{p}\pt$ Hence $\enor{S}\le1$.
\item if $p=+\infty$, for every $\ndn$,
$\eva{\ep{Sx}_n}=\eva{x_{n+1}}\le\enp{x}{\infty}$,
so $\enp{Sx}{\infty}\le\enp{x}{\infty}\pt$
\end{itemize}

Let $\alpha\in\io{0}{1}$. $x$ defined by $x_n=\alpha^n$ is in \plp{p}
and $\enp{x}{p}>0$. $Sx=\alpha x$, so $\enor{S}\ge\alpha$. $\alpha$ is as
close as wanted of $1$, so $\enor{S}=1\pt$

It remains to show that $S$ is onto. Now, if $y\in\plp{p}$, for $k\in\sC$
$x$ defined by $x_0=k$ and $x_n=y_{n-1}$ if $n\ge1$ is clearly in \plp{p},
and at once $Sf=g$. Moreover $S$ is not one-one, which proves that
$0\in\sigma\ep{S}\pt$

For $\alpha\in\sC$ such that $0<\eva{\alpha}<1$, $x$ defined by $x_n=\alpha^n$
is in \plp{p} and $x\neq0$. $Sx=\alpha x$, so $\alpha$ is an eigenvalue
of $S$, \emph{a fortiori} $\alpha\in\sigma(S)\pt$

Since $\sigma\ep{S}$ is closed, $\sigma\ep{S}\supseteq\sD$.

For $\zeta>1$, $\zeta I$ is invertible, its inverse is
$\zeta^{-1}I$ and $\enor{\ep{\zeta I}^{-1}}^{-1}=\zeta$.

Now, $\enor{S}=1<\zeta$, so, by the usual inversion theorem,
$\zeta I-S$ is invertible in \enlinco{\sE}, which proves
$\zeta\in\rho\ep{S}$ and $\zeta\notin\sigma\ep{S}\pt$
Therefore $\sigma\ep{S}=\sD$.

For $\zeta\in\rho(S)=\sC\setminus\sD$, we can apply the series expansion of
$I-\zeta^{-1}S$ which is the basis of the inversion theorem \ref{thinv-base}:
\begin{eqnarray*}
\ep{I-\zeta^{-1}S}^{-1} & = & \sum\limits^{\infty}_{k=0}\,\zeta^{-k}S^k\vrg\\
\resol(\zeta) & = & \sum\limits^{\infty}_{k=0}\,\zeta^{-k-1}S^k\pt
\end{eqnarray*}
$\resol(\zeta)=\zeta^{-1}\ep{I-\zeta^{-1}S}^{-1}\in\enlinco{\plp{p}}$,
so $\resol(\zeta)x\in\plp{p}\pt$ We can compute, for $\ndn$, since
$\ep{S^kx}_n=x_{n+k}$:

\centerline{\begin{math}\ep{\resol(\zeta)x}_n=\sum\limits^{\infty}_{k=0}%
\,\zeta^{-k-1}x_{n+k}\pt\end{math}}
\end{dem}

\begin{exemple}
\label{exemple-opcont}
With $\sE=\plp{p}$, let \sF\ the set of $x\in\sE$ such that $x_0=0$.
\sF\ is a closed vector subspace of \sE\ as kernel of the continuous linear
functional $x\,\mapsto\,x_0\pt$
Let $J$ the canonical injection from \sF\ to \sE. Let $T$ the mapping
from \sF\ to \sE\ restriction of $S$ to \sF.
\begin{enumerate}
\item $T$ is a linear isometry from \sF\ onto \sE.
\item $0\in\rho(T)\pt$
\item if $\zeta\in\sC$ and $\eva{\zeta}>1$, then $\zeta\in\sigma(T)\pt$
\end{enumerate}
\end{exemple}
\begin{dem}
$T$ is the restriction to \sF\ of $S$, hence directly
linearity and continuity of $T$ and since, for $x\in\sF$,
$\enp{Tx}{p}=\enp{Sx}{p}\le\enp{x}{p}$, $\enor{T}\le1\pt$

$T$ is onto because, among the reciprocal images of $x\in\plp{p}$ built in
the proof of prop. \ref{prop-aux-ex-opcont}, one (and only one), for $k=0$,
is in \sF: $y$ with $y_0=0$ and, for $n>1$, $y_n=x_{n-1}\pt$
Let $x\in\sF$.
\begin{itemize}
\item if $p<+\infty$, since $x_0=0$:
\begin{eqnarray*}
\enppuip{Tx}{p} & = & \sum\limits_\ndn\,\eva{x_{n+1}}^p\\
	& = & \sum\limits^{\infty}_{n=1}\,\eva{x_n}^p\\
	& = & \sum\limits^{\infty}_{n=0}\,\eva{x_n}^p-\eva{x_0}^p\\
	& = & \enppuip{x}{p}\pt
\end{eqnarray*}
So $\enp{Tx}{p}=\enp{x}{p}\pt$
\pagebreak[3.999999]
\item if $p=+\infty$, for every $\ndn$,
$\eva{\ep{Tx}_n}=\eva{x_{n+1}}\le\enp{x}{\infty}$,
so $\enp{Tx}{\infty}\le\enp{x}{\infty}\pt$

But, $\eva{x_0}=0\le\enp{Tx}{\infty}$ and
$\eva{x_n}=\eva{\ep{Tx}_{n-1}}\le\enp{Tx}{\infty}$ for $n\in\sN^\etoile$,
so $\enp{Tx}{\infty}=\enp{x}{\infty}\pt$
\end{itemize}
Hence $T$ is an isometry and $\enor{T}=1$.
Therefore $T^{-1}$ is a linear isometry from \sE\ to \sF, which implies
\refenum{ii} since $-T^{-1}=\resol(0)$.
Since \sF, closed subspace of \sE, is a Banach space endowed with
the induced  norm, $T$ is a closed operator.

Let $\zeta\in\sC$ such that $\eva{\zeta}>1$. We will search, for $y\in\sE$,
if there exists $x\in\sF$ such that $(\zeta J-T)x=y$, i.e. $\zeta x-Sx=y\pt$
Now, there exists an only solution of this equation \emph{in \plp{p}},
$x=R_{S}(\zeta)y$, verifying the equality, for each $\ndn$:
\begin{displaymath}
x_n=\sum\limits^{\infty}_{k=0}\,\zeta^{-k-1}y_{n+k}\vrg
\end{displaymath}

But have we $x\in\sF$?
Take $y$ defined by $y_0=1$ and $y_n=0$ for $n>0$ (clearly in \plp{p}).
We have $x_0=\zeta^{-1}$ and $x_n=0$ for $n>0$
But the so computed $x$ is not element of \sF, so $y$
nas no antecedent by $\zeta J-T$, which proves $\zeta\in\sigma(T)\pt$
\end{dem}

\end{document}